\documentclass{amsart}
\newcommand{\RP}{{\mathbb {RP}}}
\newcommand{\CP}{{\mathbb {CP}}}
\newcommand{\Z}{{\mathbb Z}}
\newcommand{\R}{{\mathbb R}}
\newcommand{\C}{{\mathbb C}}

\newtheorem{theorem}{Theorem}

\newtheorem{proposition}{Proposition}
\newtheorem{corollary}{Corollary}
\newtheorem{remark}{Remark}
\newtheorem{definition}{Definition}
\newtheorem{conjecture}{Conjecture}
\newtheorem{lemma}{Lemma}

\begin{document}

\author{V.A.~Vassiliev}
\address{National Research University Higher School of Economics \\
Steklov Mathematical Institute of Russian Academy of Sciences}
\email{vva@mi.ras.ru}
\thanks{Research supported by the Russian Science Foundation grant, project
16-11-10316}

\title{Local Petrovskii lacunas at parabolic singular points of wavefronts of
strictly hyperbolic PDE's}
\date{}

\begin{abstract} We enumerate the local Petrovskii lacunas (that is, the domains
of local regularity of the principal fundamental solutions of strictly
hyperbolic PDE's with constant coefficients in $\R^N$) at the {\em parabolic}
singular points of their wavefronts (that is, at the points of types $P_8^1$,
$P_8^2$, $\pm X_9$, $X_9^1$, $X_9^2$, $J_{10}^1$, $J_{10}^3$). These points
form the next difficult family of classes of the natural classification of
singular points after the so-called {\em simple} singularities $A_k, D_k, E_6,
E_7$, $E_8$, studied previously.

Also we promote a computer program counting for topologically different
morsifications of critical points of smooth functions, and hence also for local
components of the complement of a generic wavefront at its singular points.

\end{abstract}
\maketitle

\section{Introduction}

The {\em lacunas} of a hyperbolic PDE are the components of the complement of
its wavefront such that the principal fundamental solution of this equation can
be extended from any such component to a regular function in some its
neighborhood. The theory of lacunas was created by I.G.~Petrovskii
\cite{Petrovskii 45}. He has related this regularity condition to the topology
and geometry of algebraic manifolds, and gave a criterion of it in the terms of
certain homology classes of complex projective algebraic manifolds defined by
the principal symbol of the hyperbolic operator. This theory was further
developed in numerous works including \cite{Davydova 45}, \cite{Borovikov 59},
\cite{Borovikov 61}, \cite{Leray 62}, \cite{ABG 70}, \cite{ABG 73},
\cite{Gording 77}, \cite{Vassiliev 86}, \cite{Varchenko 87}, \cite{Vassiliev
92}, \cite{APLt}; for an important preceding work see \cite{Hadamard 32}. Most
of these works treat also the local aspect of the problem, explicitly
formulated in \cite{ABG 73} in the terms of {\em local lacunas} and a local
version of the Petrovskii topological condition.
\medskip

Any hyperbolic operator with constant coefficients in $\R^N$ admits a unique
fundamental solution with support in a proper cone in the half-space $\R^N_+$
of the Cauchy problem. This fundamental solution is regular (that is, locally
coincides with some smooth analytic functions) everywhere in $\R^N$ outside
some conic semialgebraic hypersurface in $\R^N_+$, called the {\em wavefront}
of our operator. We consider only {\em strictly} hyperbolic operators, which
means that the cone $A(P) \subset \check \R^N$ of zeros of the principal symbol
of our operator $P$ is non-singular outside the origin in $\check \R^N$. Here
$\check \R^N$ is the dual space of {\em momenta} with coordinates $\eta_j
\equiv \frac{1}{i}\frac{\partial}{\partial x_j},$ so that the operator $P$ is
considered as a polynomial in these variables. In this case the wavefront $W(P)
\subset \R^N_+$ is just the cone projectively dual to $A(P)$, that is, the
union of those rays from the origin in $\R_+^N$ whose orthogonal hyperplanes in
$\check \R^N$ are tangent to the cone $A(P)$. The singular points of the
wavefront (besides the origin) correspond via the projective duality to the
inflection points of $A(P)$, that is, to those points where the rank of the
second fundamental form of this cone is smaller than $N-2$. A deep
classification of these singular points was developed in the works by
V.I.~Arnold, see e.g. \cite{AVG 82}.

\begin{definition}
\label{dll} A {\it local $C^\infty$-lacuna} $($respectively, {\it holomorphic
local lacuna}$)$ at some point of the wavefront is any component of the
complement of the wavefront in a neighborhood of this point, such that the
restriction of the principal fundamental solution to this component can be
extended to a $C^\infty$-smooth function on the closure of this component
$($respectively, to an analytic function in entire neighborhood of our point).
\end{definition}

One and the same (global) component of the complement of the wavefront can to
be a local lacuna at some points of its boundary and not to be at the other
ones.

All local lacunas occurring in the neighborhoods of all singularities of
wavefronts from an initial segment of the Arnold classification (so-called
{simple} singularities) were enumerated in \cite{Vassiliev 86} and
\cite{Vassiliev 92}. In the present work we study and enumerate the holomorphic
local lacunas neighboring to the singularities whose classes the next natural
segment of this classification.

\subsection{Previous results on local lacunas}

The non-singular points of the wavefront of the operator $P$ correspond to the
points of the cone $A(P)$, at which its second fundamental form is maximally
non-degenerate. The existence and the number of local lacunas close to such
points of the wavefront can be determined in the terms of its differential
geometry, see \cite{Davydova 45} and \cite{Borovikov 59}. Namely, a component
of the complement of the wavefront at such a point is a local lacuna if and
only if the positive inertia index of the second fundamental form of the
wavefront (with the normal directed into this component) is even. A.M.~Davydova
\cite{Davydova 45} has proved the ``only if'' part of this statement: if this
signature condition is not satisfied, then already the leading term of the
asymptotics of the fundamental solution behaves as a half-integer (but not
integer) power of the distance from the wavefront. V.A.~Borovikov
\cite{Borovikov 59}, using complicated analytic estimates, has proved that
otherwise we have a local lacuna, that is, {\em all} terms of the asymptotic
expansion of this solution in the terms of this distance have integer powers,
and the corresponding power series does converge. His result was later
explained in \cite{ABG 73} as a corollary of the removable singularity theorem
by moving into the complex domain.

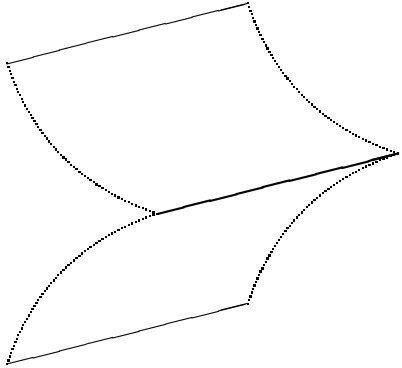
\begin{figure}
\begin{center}
\unitlength 1.00mm \linethickness{0.4pt}
\begin{picture}(50,50)
\bezier{60}(0,0)(5,15)(20,20) \bezier{60}(0,40)(5,25)(20,20)
\put(0,0){\line(4,1){32}} \put(0,40){\line(4,1){32}}
\bezier{60}(32,8)(37,23)(52,28) \bezier{60}(32,48)(37,33)(52,28) {\thicklines
\put(20,20){\line(4,1){32}}}
\end{picture}
 \caption{Cuspidal edge in the 3-dimensional space}
\label{a2}
\end{center}
\end{figure}

All local lacunas neighboring to the simplest singular points of the
wavefronts, of types $A_2$ (cuspidal edges, see Fig.~\ref{a2}) and $A_3$
(swallowtails) were counted for in \cite{Gording 77}. An interesting situation
occurs close to a point of the cuspidal edge, if $N$ is odd and the inertia
indices of the quadratic part of the {\em generating function} of our point
(see \S \ref{gengen} below) also is odd. For all other combinations of these
numbers, if a component of the complement of the wavefront close to the
cuspidal edge is not a local lacuna, then already the Davydova--Borovikov
signature condition from the side of this component is not satisfied at some
non-singular points of the wavefront arbitrarily close to the edge. However, in
the case of odd $N$ and $i_{\pm}$ the Davydova--Borovikov condition from the
side of the bigger component (see Fig.~\ref{a2}) is satisfied at all nearby
non-singular points, nevertheless this component is not a local lacuna (and
also is not for all other combinations of $N$ and $i_{\pm}$).

\begin{table}
\begin{center}
\caption{Numbers of local lacunas at simple singularities of wavefronts}
\label{t12}
\begin{tabular}{|l|c|c|c|c|}
\hline Singularity & $N$ even & $N$ even & $N$ odd & $N$ odd \cr class & $i_+$
even & $i_+$ odd & $i_+$ even & $i_+$ odd \cr \hline $A_1$ & 2 & 0 & 1 & 1 \cr
$A_{2k}, \; k\ge 1$ & 0 & 0 & 1 & 0 \cr $\pm A_{2k+1}, \; k \ge 1$ & 0 & 1 & 1
& 1 \cr $D_4^-$ & 0 & 3 & 1 & 1 \cr $D_{2k}^+, \; k \ge 2$ & 0 & 0 & 1 & 1 \cr
$D_{2k}^-, \; k \ge 3$ & 0 & 2 & 1 & 1 \cr $\pm D_{2k+1}, \; k \ge 2$ & 0 & 0 &
1 & 1 \cr $\pm E_6$ & 0 & 0 & 1 & 1 \cr $E_7$ & 0 & 0 & 1 & 1 \cr $E_8$ & 0 & 0
& 1 & 1 \cr \hline
\end{tabular}
\end{center}
\end{table}

All local lacunas for all {\em simple} singularities of wavefronts (that is,
singularities of classes $A_k, D_k,$ $E_6,$ $E_7,$ $E_8$ in the Arnold's
classification) were found in \cite{Vassiliev 86}, see Table \ref{t12} for the
number of them.

Atiyah, Bott and G\aa rding \cite{ABG 73} have introduced the local version of
the homological Petrovskii criterion, and proved that it implies that the
corresponding local component of the complement of the wavefront is a
holomorphic local lacuna. In \cite{Vassiliev 86} the converse implication was
proved for {\em finite type} points of wavefronts (that is, for points
corresponding by the projective duality to only finitely many lines in the
complexification of $A(P)$; this condition is satisfied for all singular points
of wavefronts of {\em generic} operators). In \cite{Vassiliev 92}, an easy
geometric criterion for a component to be a local lacuna of a simple
singularity was proved. Namely, it follows from the above described facts that
if a local component of the complement of the wavefront is a local lacuna, then
the Davydova--Borovikov signature condition is satisfied at all non-singular
points of its boundary, and in addition our component is the ``smaller''
component of the complement of the wavefront at all points of type $A_2$ (that
is, cuspidal edges) of this boundary, see Fig.~\ref{a2}. In \cite{Vassiliev 92}
it was proved that if the singularity is simple, and some technical condition
(the {\em versality} of the generating family, which always holds for
wavefronts of generic operators) is satisfied, then this necessary condition is
also sufficient; moreover, in this case the notions of local $C^\infty$-lacunas
and holomorphic local lacunas are equivalent.

\begin{table}
\begin{center}
\caption{Numbers of local lacunas at parabolic singularities} \label{t13}
\begin{tabular}{|l|c|c|c|c|}
\hline Singularity & $N$ even & $N$ even & $N$ odd & $N$ odd \cr class & $i_+$
even & $i_+$ odd & $i_+$ even & $i_+$ odd \cr \hline $P_8^1$ & $0_c$ & $0_c$ &
$\ge 2$ & $0$ \cr \hline $P_8^2$ & $0_c$ & $0_c$ & $\ge 2$ & 0 \cr \hline $\pm
X_9$ & $1_c$ & 0 & $\ge 2$ & 0 \cr \hline $X_9^1$ & 0 & 0 & 0 & 0 \cr \hline
$X_9^2$ & 0 & $\ge 4$ & 0 & 0 \cr \hline $J_{10}^3$ & $0_c$ & $\ge 1$ & 0 & 0
\cr \hline $J_{10}^1$ & $0_c$ & $0_c$ & 0 & 0 \cr \hline
\end{tabular}
\end{center}
\end{table}

The next important natural set of singularity classes of wavefronts is that of
{\it parabolic} (or {\it simple-elliptic}) singularities, see \cite{AVG 82}. It
consists of seven one-parameter families of singularities listed in the
left-hand column of Table \ref{t13}. The numbers of local lacunas at these
singularities are shown in the remaining columns of this table: this (together
with an explicit description of these lacunas) is the main result of the
present article, see Theorem \ref{thm1} below. However, we need some
preliminaries to describe these singularities and formulate this result
accurately.

\subsection{Generating functions and generating families of wavefronts}
\label{gengen}

Given a point $x \in \R^N \setminus 0$ of the wavefront of a strictly
hyperbolic operator $P$, the local geometry of this wavefront at this point (in
particular the set of local components of its complement at this point) is
determined by its {\em generating function}, which is just the function $f$ in
the local equation
$$\xi_{0} = f(\xi_1, \dots, \xi_{N-2}) $$ of the projectivization $A^*(P)
\subset \check{\RP}^{N-1}$ of the set of zeros of the principal symbol of our
operator. Here $\xi_0, \dots, \xi_{N-2}$ are affine local coordinates in
$\check \RP^{N-1}$ with the origin at the tangency point of the hypersurface
$A^*(P)$ and the hyperplane $L(x)$ orthogonal to the line containing the point
$x$, such that this hyperplane $L(x)$ is distinguished by the equation
$\xi_0=0$. In particular, $f$ has a critical point at the origin; the dual
piece of the wavefront is smooth if this critical point is Morse.

Denote by $n$ the number $N-2$ of variables of generating functions of
wavefronts in $\R^N$. The {\em parabolic} singularity classes studied in this
work have the generating functions which can be reduced by a local
diffeomorphism in $\R^{n}$ to the following normal forms. The functions of
class $P_8$ in appropriate (curvilinear) local coordinates have the formula
$\varphi(x_1, x_2, x_3)+Q(x_4, \dots, x_n),$ where $\varphi$ is a
non-degenerate homogeneous cubic polynomial, and $Q$ is a non-degenerate
quadratic function in the remaining coordinates, e.g. $\pm x_4^2 \pm \dots \pm
x_n^2$. The projectivization of the zero set of the polynomial $\varphi$ can
consist of one or two curves, therefore we obtain two subclasses, called
$P_8^1$ and $P_8^2$ respectively. The remaining parabolic functions have the
following normal forms (where $Q$ are non-degenerate quadratic functions in
coordinates $x_3, \dots, x_n$):

$$
\begin{array}{cll}
\pm X_9 & \pm (x_1^4 + \alpha x_1^2 x_2^2 + x_2^4 + Q) & \alpha > -2 \\
X_9^1 &  x_1x_2(x_1^2 + \alpha x_1x_2 + x_2^2) + Q & \alpha^2 < 4  \\
X_9^2 &  x_1x_2(x_1 + x_2)(x_1 + \alpha x_2) + Q & \alpha \in (0,1) \\
J_{10}^3 & x_1(x_1 - x_2^2)(x_1 - \alpha x_2^2) + Q & \alpha \in (0,1) \\
J_{10}^1 & x_1(x_1^2 + \alpha x_1x_2^2 + x_2^4) + Q & \alpha^2 < 4
\end{array}
$$

The index $i_+$ in Table \ref{t13} is the positive inertia index of the
quadratic part $Q$ of the corresponding function.

Another important notion, reducing the study of wavefronts to the context of
critical points of functions, is that of {\em generating families}. In our
case, this is the name of the family of functions \begin{equation}
\label{genfam} f_\lambda \equiv f(\xi_1, \dots, \xi_{N-2}) -\lambda_0 -
\lambda_1 \xi_1 - \dots - \lambda_{N-2}\xi_{N-2},\end{equation} depending on
the parameter $\lambda =(\lambda_0, \dots, \lambda_{N-2}) \in \R^{N-1}$. It is
natural to consider these parameters $\lambda_i$ as local affine coordinates in
$\RP^{N-1}$ close to the point $\{x\}$. Indeed, any collection $\lambda$ of
these numbers defines a hyperplane $L(\lambda) \subset \check \RP^{N-1}$
distinguished by the equation $$\xi_0=\lambda_0 + \lambda_1 \xi_1 + \dots +
\lambda_{N-2}\xi_{N-2},$$ hence a line in $\R^N$ or a point in $\RP^{N-1}$. The
projectivized wavefront close to our point in $\RP^{N-1}$ consists of all {\em
discriminant} values of the parameters of the family (\ref{genfam}), that is,
of those values of $\lambda$ for which the function (\ref{genfam}) has critical
value 0 (which is equivalent to the tangency of hypersurfaces $A^*(P)$ and
$L(\lambda)$). So, the role of the (projectivized) wavefronts in the language
of critical points of functions is played by the {\em discriminant varieties}
of {\em function deformations}.

Recall that a {\em deformation} of the function $f: \R^n \to \R$ is a function
$F: \R^n \times \R^l \to \R$ considered as a family of functions $f_\lambda
\equiv F(\cdot ,\lambda): \R^n \to \R$ depending on the parameter $\lambda \in
\R^l$, such that $f_0$ coincides with the deformed function $f$. The
discriminant variety of such a deformation is the set of values of its
parameter $\lambda$ such that the corresponding function $f_\lambda$ has a
critical point with zero critical value. In particular, the generating family
of a wavefront is a deformation of its generating function, and the wavefront
itself is the discriminant set of this deformation.

The notion of (holomorphic) local lacunas has sense for arbitrary deformations
of critical points of real functions (not necessarily related with the
wavefronts): a component of the complement of the discriminant set of such a
deformation is a local lacuna if some homological condition (the triviality of
the local Petrovskii homology class, described in the next section) concerning
the level manifolds $f_\lambda^{-1}(0)$ is satisfied for values of $\lambda$
from this component. If our deformation is the generating family of the
wavefront of a hyperbolic operator, then this notion turns out to be equivalent
to the one described previously in the terms of fundamental solutions, see
Proposition \ref{petcrit} of the next section.

Among the deformations $F(x,\lambda)$ of a function $f$ with an isolated
critical point there is a distinguished class of {\em versal deformations},
that is, of sufficiently ample deformations such that all other deformations
can be reduced to them in some precise sense, see \cite{AVG 82}, \cite{APLt}.
The number of parameters of a versal deformation cannot be smaller than the
{\em Milnor number} $\mu(f)$ of our critical point (which is the standard lower
index of the notation of its class, like 8 for $P_8^1$), on the other hand
almost all deformations depending on $\ge \mu(f)$ parameters satisfy this
condition. An important corollary of the notion of versality is as follows: if
the parameter space of some versal deformation of the function $f$ contains no
local lacunas, then the same is true for any other its deformation.

\begin{theorem}\label{thm1}
The number of holomorphic local lacunas in the parameter space of a versal
deformation of a parabolic critical point of a function $f_0:\R^{N-2} \to \R$
is equal to the number indicated in the corresponding cell of Table \ref{t13}
or satisfies the inequality given in this table. $($The subscript $ _c$ in some
cells means that in the corresponding case the upper bound on the number of
local lacunas has only a computer proof$)$.

In the case of non-versal deformations this statement remains true for all
cells of the Table \ref{t13}, where we have $0$ or $0_c$. If our deformation
contains all functions $f_0 + \mbox{const}$ $($which holds for all generating
families of wavefronts$)$ then both signs $1_c$ and $\ge 2$ for the singularity
$\pm X_9$ can be replaced by $\ge 1$.
\end{theorem}

\begin{remark} \rm Some of these results were obtained earlier, see
e.g. \cite{AVGL 89}, \cite{APLt}. The new results are as follows:

The singularity $P_8^2$ is investigated for the first time.

For $\pm X_9$, $N$ even, $i_+$ even: the estimate $ \ge 1$ is replaced by the
exact value $1_c$.

For $X_9^2$, $N$ even, $i_+$ odd: $ \ge 2$ is replaced by $\ge 4$.

For $J_{10}^3$, $N$ even, $i_+$ even: the absence of local lacunas is proved.

For $J_{10}^1$, $N$ even, $i_+$ even or odd: the absence of local lacunas is
proved in both cases.
\end{remark}

All local lacunas mentioned in non-zero cells of Table \ref{t13} will be
presented in \S \ref{realiz}. All zeros (but not signs $0_c$) in this table
follow from a topological obstruction described in \S \ref{bound}. All signs
$0_c$ are proved by a combinatorial Fortran program which enumerates all
possible topological types of morsifications of given critical points and
checks the local Petrovskii condition for them. This program has also found for
the first time the local lacunas for singularities $P_8^1$, $X_9^2$ and
$J_{10}^3$ presented below, as well as one of local lacunas for $\pm X_9$, see
Proposition \ref{otherx9} on the page \pageref{otherx9}. This program is
described in \S \ref{progr}. The upper bound $1_c$ for the singularity $\pm
X_9$ with even $N$ and $i_+$ will be proved in \S \ref{noextra}.

\begin{conjecture}
In all cells of Table \ref{t13} $($except maybe for the case $X_9^2)$ the
inequalities can be replaced by equalities.
\end{conjecture}

Here is a particular information supporting this conjecture. Given a Morse
perturbation $f_\lambda$ of a function $f$ with complicated critical point at
the origin, denote by $\chi(\lambda)$ the number of its real critical points
with negative critical values and even Morse index minus the number of points
also with negative values but odd Morse index.

\begin{proposition} \label{chi} The perturbations of our singularity $f$, which belong to
the local lacunas, cannot have values of $\chi(\lambda)$ different from those
for perturbations presented in \S \ref{realiz}.
\end{proposition}

This fact is also proved by our program. \hfill $\Box$

\section{Local Petrovskii classes and their properties}

Let $f: (\C^n, \R^n, 0) \to (\C,\R,0)$ be a holomorphic function with isolated
critical point at $0$, $\mu(f)$ its {\em Milnor number} (see \cite{AVG 82}),
$B_\varepsilon \subset \C^n$ a ball centered at $0$ with a small radius
$\varepsilon$. Let $f_\lambda$ be a very small (with respect to $\varepsilon$)
perturbation of $f$, such that $0$ is not a critical value of $f_\lambda$ in
$B_\varepsilon$. Consider the corresponding {\it Milnor fiber} $V_\lambda
\equiv f_\lambda^{-1}(0) \cap B_\varepsilon$. It is a smooth
$(2n-2)$-dimensional manifold with boundary $\partial V_\lambda \equiv
V_\lambda \cap \partial B_\varepsilon$. By the Milnor's theorem it is homotopy
equivalent to the wedge of $\mu(f)$ spheres $S^{n-1}$, in particular $\tilde
H_{n-1}(V_\lambda) \simeq \Z^{\mu(f)} \simeq \tilde H_{n-1}(V_\lambda,
\partial V_\lambda)$; here $\tilde H_*(V_\lambda)$ denotes the homology group
reduced modulo a point, and $\tilde H_*(V_\lambda, \partial V_\lambda)$ the
relative homology group reduced additionally modulo the fundamental cycle.
Also, $\mu(f)$ is equal to the number of critical points of $f_\lambda$ in
$B_\varepsilon$, if the function $f_\lambda$ is {\em Morse} and is indeed
sufficiently close to $f$.

We will assume that an orientation of $\R^n$ is fixed, and the differential
form $dx_1 \wedge \dots \wedge dx_n$ is positive with respect to this
orientation.

\begin{figure}\unitlength 0.850mm
\linethickness{0.4pt}
\begin{picture}(143.00,62)
\put(43.00,31.00){\circle*{1.33}} \put(26.50,31.00){\circle*{1.33}}
\put(12.00,31.00){\circle*{1.33}} \put(71.00,31.00){\makebox(0,0)[cc]{\small
$\R^1$}} \put(49.00,41.00){\circle*{1.33}} \put(49.00,21.00){\circle*{1.33}}
\put(20.00,16.00){\circle*{1.33}} \put(20.00,46.00){\circle*{1.33}}
\put(0.00,32.00){\vector(1,0){10.00}} \put(14.00,32.00){\vector(3,-4){0.20}}
\bezier{36}(10.00,32.00)(12.00,36.00)(14.00,32.00)
\put(14.00,32.00){\vector(1,0){9.00}}
\bezier{56}(23.00,32.00)(26.50,25.00)(30.00,32.00)
\put(30.00,32.00){\vector(1,0){11.00}} \put(45.00,32.00){\vector(3,-4){0.20}}
\bezier{36}(41.00,32.00)(43.00,36.00)(45.00,32.00)
\put(45.00,32.00){\vector(1,0){17.00}} \put(0.00,30.00){\vector(1,0){10.00}}
\put(14.00,30.00){\vector(3,4){0.20}}
\bezier{36}(10.00,30.00)(12.00,26.00)(14.00,30.00)
\put(14.00,30.00){\vector(1,0){9.00}}
\bezier{72}(23.00,30.00)(26.50,37.00)(30.00,30.00)
\put(30.00,30.00){\vector(1,0){11.00}} \put(45.00,30.00){\vector(3,4){0.20}}
\bezier{36}(41.00,30.00)(43.00,26.00)(45.00,30.00)
\put(45.00,30.00){\vector(1,0){17.00}} \put(0.00,31.00){\vector(1,0){63.00}}
\put(26.50,28.50){\vector(1,0){1.00}} \put(26.50,33.50){\vector(1,0){1.00}}

\put(123.00,31.00){\circle*{1.33}} \put(106.50,31.00){\circle*{1.33}}
\put(92.00,31.00){\circle*{1.33}}  \put(129.00,41.00){\circle*{1.33}}
\put(129.00,21.00){\circle*{1.33}} \put(100.00,16.00){\circle*{1.33}}
\put(100.00,46.00){\circle*{1.33}} \put(80.00,31.00){\vector(1,0){63.00}}

\put(100,13){\vector(-1,0){1.00}} \put(100,49){\vector(-1,0){1.00}}
\put(129,18){\vector(-1,0){1.00}} \put(129,44){\vector(-1,0){1.00}}
\put(100,16){\circle{6}} \put(100,46){\circle{6}} \put(129,41){\circle{6}}
\put(129,21){\circle{6}} 

\bezier{90}(0,31)(0,0)(31.50,0) \bezier{90}(31.50,0)(63,0)(63,31)
\bezier{90}(63,31)(63,62)(31.50,62) \bezier{90}(31.50,62)(0,62)(0,31)
\put(33,57){\small $B_\varepsilon$}

\bezier{90}(80,31)(80,0)(111.50,0) \bezier{90}(111.50,0)(143,0)(143,31)
\bezier{90}(143,31)(143,62)(111.50,62) \bezier{90}(111.50,62)(80,62)(80,31)
\put(113,57){\small $B_\varepsilon$}

\end{picture}
\caption{Odd Petrovskii cycle for $n=1$} \label{oddp}
\end{figure}
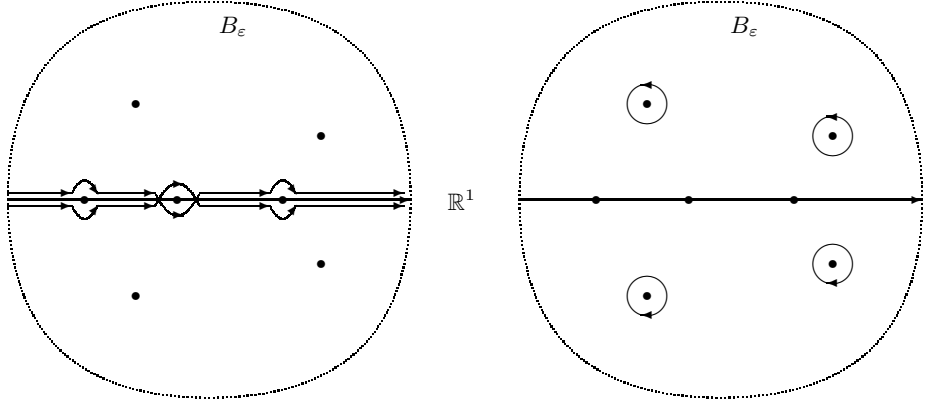

There are two important elements in the group $\tilde H_{n-1}(V_\lambda,
\partial V_\lambda)$, the {\em even} and {\em odd Petrovskii classes}. The
first of them, $P_{\mbox{ev}}(\lambda),$ is presented by the cycle of real
points $\R^n \cap V_\lambda$ oriented by the differential form $(dx_1 \wedge
\dots \wedge dx_n) /df_\lambda$. The definition of the second class,
$P_{\mbox{odd}}$, is a bit more complicated. First we consider an
$n$-dimensional cycle $\Pi(\lambda)$ in $B_\varepsilon \setminus V_\lambda$,
presented by two copies of canonically oriented $\R^n$, slightly moved in a
small neighborhood of the submanifold $\R^n \cap V_\lambda$ in $B_\varepsilon$
so that they streamline $V_\lambda$ from two different sides in the complex
domain: for the case $n=1$ see the left-hand part of Fig. \ref{oddp} where the
set $V_\lambda$ is marked by thick dots.

The {\em odd Petrovskii class} $P_{\mbox{odd}}$ is defined as the preimage of
the homology class of this cycle under the {\it Leray tube operator} $$\tilde
H_{n-1}(V_\lambda,
\partial V_\lambda) \to \tilde H_n(B_\varepsilon \setminus V_\lambda, \partial
B_\varepsilon),$$ which sends any relative cycle in the submanifold $V_\lambda$
to the union of boundaries of the fibers of the tubular neighborhood of this
submanifold over the points of this cycle. This operator is conjugate via the
Poincar\'e--Lefschetz isomorphisms to the boundary isomorphism $\tilde
H_n(B_\varepsilon, V_\lambda) \to \tilde H_{n-1}(V_\lambda),$ and also is an
isomorphism.
\medskip

Let $F:(\R^n \times \R^l,0) \to (\R,0)$ be a deformation of the function $f$,
and $\Sigma(F) \subset \R^l$ be the set of discriminant values of the parameter
$\lambda$.

\begin{definition} \rm
A local (close to the point $0 \in \R^l$) connected component of the set $\R^l
\setminus \Sigma(F)$ is called an {\em even} (respectively, {\em odd}) {\em
local lacuna} of the deformation $F$ if for any value of $\lambda$ from this
component the element $P_{\mbox{ev}}(\lambda)$ (respectively,
$P_{\mbox{odd}}(\lambda)$) of the group $\tilde H_{n-1}(V_\lambda, \partial
V_\lambda)$, related with the corresponding perturbation $f_\lambda$ of $f$, is
equal to 0.
\end{definition}

This definition is consistent with Definition \ref{dll} by the following
reason.

\begin{proposition} \label{petcrit} Suppose that $x \in \R^N \setminus 0$ is a
point of the wavefront of a strictly hyperbolic operator, the critical point of
its generating function $f$ is isolated, and $y$ is a point outside the
wavefront but very close to $x$. Let $y^* \in \RP^{N-1}$ be the direction of
the line containing $y$, $\lambda=(\lambda_0, \dots, \lambda_{N-2})$ the local
coordinates of the point $y^*$ accordingly to \S \ref{gengen}, and $f_\lambda$
the perturbation of $f$ given by the formula $($\ref{genfam}$)$ with these
values of $\lambda_i$. Then the point $y$ belongs to a local $($close to $x)$
holomorphic lacuna of our hyperbolic operator if and only if the point
$\lambda$ belongs to an even even lacuna of the corresponding generating family
$($if $N$ is even$)$ or to an odd local lacuna $($if $N$ is odd$)$.
\end{proposition}

The part ``if'' of this proposition was essentially proved in \cite{ABG 73},
and ``only if'' was conjectured there and proved in the translator's note to
the Russian translation of \cite{ABG 73}, see also \cite{Vassiliev 86}.

This proposition reduces the study of local lacunas to a problem on
deformations of real critical points of functions, namely to the calculation of
local Petrovskii classes of perturbations of such functions, and the hunt for
those perturbations for which these classes vanish.

\subsection{Important example} \label{examin} If the function $f$ has a minimum
$($res\-pec\-ti\-vely,  maximum$)$ point at $0$, then the function $f+\tau$
$($respectively, $f-\tau)$ with sufficiently small positive $\tau$ belongs to
its even local lacuna.
\medskip

Indeed, in this case the set of real points of the corresponding Milnor fiber
is empty.

\begin{conjecture} If a function $f(x_1, x_2)$ has an isolated
non-Morse critical point at 0, then its deformations have no even local lacunas
unless $f$ has an extremum at the origin; in the latter case all critical
values of real critical points of all its small perturbations which belong to
such a lacuna are positive $($if $f$ has a minimum point at the origin$)$ or
negative $($if $f$ has a maximum$)$.
\end{conjecture}

\subsection{Explicit calculation of local Petrovskii classes}
\label{expl}

The group $\tilde H_{n-1}(V_\lambda, \partial V_\lambda)$ is Poincar\'e dual to
$\tilde H_{n-1}(V_\lambda)$, therefore any its element is completely
characterized by its intersection indices with basic elements of the latter
group. For these basic elements we can take the {\em vanishing cycles} (see
e.g. \cite{AVG 84}, \cite{Gusein-Zade 77}, \cite{APLt}) corresponding to
critical points of $f_\lambda$. In \cite{Leray 62}, \cite{Vassiliev 86}
explicit formulas for these intersection indices of both local Petrovskii
classes with cycles vanishing in real critical points were calculated: these
indices are expressed in the terms of Morse indices of these critical points
and the intersection indices of vanishing cycles. We do not give here these
large formulas and refer to section V.1.6 in \cite{APLt} or \S 5.1.4 in
\cite{AVGL 89}.

In particular, these formulas describe the Petrovskii classes completely if all
$\mu(f)$ critical points of the perturbation $f_\lambda$ are real, and we know
their Morse indices and the intersection indices of corresponding vanishing
cycles.

\begin{remark} \rm
In this and other related calculations it is important to use the orientations
of these vanishing cycles, compatible with the fixed orientation of $\R^n$.
Fortunately, the methods of calculating the intersection indices developed in
\cite{Gusein-Zade 74}, \cite{A'Campo 75}, \cite{Ga1} use exactly these
orientations. These methods give us the desired data for appropriate
perturbations of all parabolic singularities except for $P_8^2$, in the latter
case we solve the similar problem in \S \ref{DDP82}.
\end{remark}

\subsection{Stabilization}

\begin{definition}[see \cite{AVG 82}] \rm Two functions $f, \tilde f:(\R^n,0)
\to (\R,0)$ with critical points at $0$ are {\em equivalent} if they can be
taken one into the other by a germ of diffeomorphism $G:(\R^n,0)\to (\R^n,0)$,
that is, $f \equiv \tilde f\circ G.$ Two critical points of functions (maybe
depending on a different number of variables) are {\em stably equivalent} if
they become equivalent after the summation with non-degenerate quadratic forms
depending on additional variables.
\end{definition}

For example, the functions $f(x)=x^k$, $f_1(x,y)=x^k +y^2$, $f_2(x,y,z)=x^k
-yz$ and $f_3(x,y)=-y^2+(x-y)^k$ are stably equivalent to one another, but they
are not stably equivalent to the function $f_4(x,y)=x^k$.

If $F(x,\lambda)$ is a deformation of the function $f(x)$, $x=(x_1, \dots,
x_n)$, then the family of functions $F(x_1, \dots, x_n,\lambda)\pm x_{n+1}^2
\pm \dots \pm x^2_{n+m}$ depending on $n+m$ variables is a deformation of the
stabilization $f(x) \pm x^2_{n+1} \pm \dots \pm x^2_{n+m}$ of $f(x)$; the
latter deformation is versal if and only if $F(x,\lambda)$ is.

\begin{proposition}[see e.g. \cite{APLt}] \label{stabil} The following
conditions are equivalent:

1$)$ the perturbation $f_\lambda$ of the function $f:(\C^n,\R^n,0)\to(\C,\R,0)$
belongs to an even $($respectively, odd$)$ local lacuna;

2$)$ the perturbation $f_\lambda+x_{n+1}^2 +x_{n+2}^2$ of the function
$f+x_{n+1}^2+x_{n+2}^2:(\C^{n+2},\R^{n+2},0)\to(\C,\R,0)$ belongs to an even
$($respectively, odd$)$ local lacuna;

3$)$ the perturbation $f_\lambda-x_{n+1}^2 -x_{n+2}^2$ of the function
$f-x_{n+1}^2-x_{n+2}^2:(\C^{n+2},\R^{n+2},0)\to(\C,\R,0)$ belongs to an even
$($respectively, odd$)$ local lacuna;

4$)$ the perturbation $f_\lambda+x_{n+1}^2 -x_{n+2}^2$ of the function
$f+x_{n+1}^2-x_{n+2}^2:(\C^{n+2},\R^{n+2},0)\to(\C,\R,0)$ belongs to an
\underline{odd} $($respectively, \underline{even}$)$ local lacuna.
\end{proposition}

So, the summation with a positive or negative definite quadratic form in an
even number of additional variables moves even (respectively, odd) local
lacunas to the lacunas of the same type; the summation with a quadratic
function of signature $(1,1)$ in two additional variables moves even lacunas to
odd ones and vice versa. Therefore any stable equivalence class of functions
splits into four subclasses, depending on the parities of $n$ and of an
arbitrary (say, positive) inertia index of the quadratic part of the Taylor
expansion of these functions. The sets of local lacunas (of the same parity) of
versal deformations of functions from any of these four subclasses are in a one
to one correspondence with each other.

\begin{corollary}
\label{extra} If the function $f(x_1, \dots, x_n)$ has a minimum
$($respectively, maximum$)$ point at $0$, then the function $f(x_1, \dots, x_n)
-x_{n+1}^2 +\tau$ $($res\-pec\-ti\-vely, $f(x_1, \dots, x_n) + x_{n+1}^2
-\tau)$ with sufficiently small $\tau>0$ belongs to an {\em odd} local lacuna.
\end{corollary}

Indeed, say in the first case the function $f -x_{n+1}^2+x_{n+2}^2-x_{n+3}^2
+\tau$ belongs to an even local lacuna by item 3) of Proposition \ref{stabil}
and by \S \ref{examin}; it remains to use item 4) of the same proposition.
\hfill $\Box$

\subsection{Multiplication by $-1$}
\label{pm}

It follows immediately from the definitions of Petrovskii classes that the
perturbation $-f_\lambda$ of the function $-f$ belongs to a local lacuna if and
only if the perturbation $f_\lambda$ of $f$ does.

\subsection{Another form of the odd Petrovskii cycle} \label{another}

It is easy to see that the cycle of Fig.~\ref{oddp} (left) is homological in
$\C^1 \setminus V_\lambda$ (modulo the complement of $B_\varepsilon$) to the
sum of small circles going around all non-real points of the set $V_\lambda$,
as shown in the right-hand part of this figure. So the pre-image of its
homology class under the Leray tube operator, that is, the odd Petrovskii
class, is represented by the sum of these points taken with appropriate signs.
The same construction allows us to realize this class in the case of an
arbitrary odd $n$.

Namely, let us choose a point ${\bf x} \in \R^n \setminus V_\lambda$. Let
$S({\bf x}) \sim S^{n-1}$ be the space of all oriented affine lines in $\R^n$
through ${\bf x}$, $\phi: E({\bf x}) \to S({\bf x})$ the tautological line
bundle, and $\phi_{\C}: E_{\C}({\bf x}) \to S({\bf x})$ its complexification,
so that $E_\C$ is the union of pairs $(l,x) \in S({\bf x}) \times \C^n$ such
that $x$ belongs to the complexification of the line $\{l\} \subset \R^n$. The
manifold $E$ is obviously orientable and is separated into two parts by the
section of the line bundle consisting of all points $(l,{\bf x})$. The
forgetful map $\Psi:(l,x) \mapsto x$ sends any of these parts diffeomorphically
to $\R^n \setminus {\bf x}$. In the case of odd $n$ the orientations of these
parts induced from the fixed orientation of $\R^n$ belong to one and the same
orientation of entire $E$.\footnote{In the case $n=1$, responsible for
Fig.~\ref{oddp}, the role of the orientation of the base is played by the
choice of (different) signs of two points of the 0-dimensional sphere $S({\bf
x})$. So the canonical orientation of the line over the negative point should
be reversed.} Extend the map $\Psi$ to the similar forgetful map $\Psi_\C :
E_\C \to \C^n$. For any oriented line $l \in S({\bf x})$ the set
$\Psi_{\C}^{-1}(V_\lambda) \cap \{l_\C\}$ is a finite set symmetric with
respect to the real line $\{l\} \subset \{l_\C\}$. Define the relative cycle
$\tilde P({\bf x}) \subset E_{\C} \cap \Psi_{\C}^{-1}(B_\epsilon \setminus
V_\lambda) $ as the union (over all points $l \in S({\bf x})$) of real lines
$\{l\}$ slightly moved inside their complexifications $\{l_\C\}$ close to all
points of $\Psi^{-1}(V_\lambda)$ in such a way that they bypass these points
from the left side with respect to the canonical orientation of $\{l\}$. The
homology class of the cycle $\Pi(\lambda)$ from the construction of the odd
Petrovskii class can be realized as the direct image of this cycle $\tilde
P({\bf x})$ under the map $\Psi_{\C}$. (In fact, $\tilde P({\bf x})$ is a kind
of ``blowing up'' the cycle $\Pi(\lambda)$ at the point ${\bf x}$).

For any $l$ the obtained 1-dimensional cycle in $\{l_\C\} \setminus
\Psi^{-1}(V_\lambda)$ is homological (within the upper half-plane and modulo
the intersection with $\Psi^{-1}(\C^n \setminus B_\varepsilon)$) to the union
of small circles around all imaginary points of $\Psi^{-1}(V_\lambda \cap
B_\varepsilon)$ in this left-hand half-plane. These homologies can be performed
uniformly over all $l$ and sweep out a homology between the cycle $\tilde
P({\bf x})$ and a cycle which is the Leray tube around the union (over all $l
\in S({\bf x})$) of all such imaginary points with positive imaginary parts in
the fibers $\{l_\C\}$. Thus the odd Petrovskii class can be realized in the
case of odd $n$ as the direct image under the map $\Psi_\C$ of the cycle
composed by this union.

This is essentially the original definition of the odd Petrovskii cycle, see
\cite{Petrovskii 45}. Unlike the even cycle, it depends on the choice of the
point ${\bf x}$, but its homology class does not.

\section{An obstruction to the existence of local lacunas}
\label{bound}

Either of two local Petrovskii classes related to a non-discriminant point
$\lambda \in \R^l$ is an element of the relative homology group $\tilde
H_{n-1}(V_\lambda, \partial V_\lambda)$ of the corresponding Milnor fiber
$V_\lambda = f_\lambda^{-1}(0) \cap \partial B_\epsilon$. The boundary operator
of the exact sequence of the pair $(V_\lambda, \partial V_\lambda)$ sends this
class to some element of the group $\tilde H_{n-2}(\partial V_\lambda)$. For
any deformation $F(x,\lambda)$ of the function $f(x)$, the spaces $\partial
V_\lambda$ form a locally trivial (and hence trivializable) fiber bundle over a
neighborhood of the origin in the parameter space of our deformation (including
the discriminant values of $\lambda$). Therefore the homology groups $\tilde
H_{n-2}(\partial V_\lambda)$ over all values of $\lambda$ are naturally
identified to one another. It follows easily from the construction of the
Petrovskii classes, that the boundaries of all cycles $P_{\mbox{ev}}(\lambda)$
(respectively, $P_{\mbox{odd}}(\lambda)$) over all non-discriminant values of
$\lambda$ are mapped into one another by this identification. Therefore if for
some value $\lambda \in \R^l$ this boundary is not homologous to zero, then the
same is true for all other values of $\lambda$, in particular the local
Petrovskii classes of the same parity for all $\lambda$ are non-trivial and we
have no local lacunas for the corresponding singularity.

All zeros in the cells of Table \ref{t13} (but not the signs $0_c$) except for
the last column of $P_8^2$ follow from this obstruction and from the explicit
calculation of the Petrovskii classes mentioned in \S \ref{expl}. Conversely,
in the case of $P_8^2$ the calculation of intersection indices of vanishing
cycles in \S \ref{DDP82} will be based on the calculation of this boundary
which will be done in \S \ref{obstrp8}.

\begin{remark} \rm
\label{rem1} The group $\tilde H_*(\partial V_\lambda)$ can be non-trivial only
in dimensions $n-1$ and $n-2$, and its structure can be obtained from the
intersection form in $\tilde H_{n-1}(V_\lambda)$. Indeed, by the Milnor's
theorem the only non-trivial segment of the exact sequence of the pair
$(V_\lambda,
\partial V_\lambda)$ is
\begin{equation}
\label{exact} 0 \to \tilde H_{n-1}(\partial V_\lambda) \to \tilde
H_{n-1}(V_\lambda) \stackrel{j}{\to} \tilde H_{n-1}(V_\lambda, \partial
V_\lambda) \to \tilde H_{n-2}(\partial V_\lambda) \to 0.
\end{equation}
If we fix Poincar\'e dual frames in two central groups of this sequence (which
are isomorphic to $\Z^{\mu(f)}$), then the homomorphism $j$ will be given by
the intersection matrix of basic elements of $\tilde H_{n-1}(V_\lambda)$. It
determines completely both marginal groups $\tilde H_i(\partial V_\lambda)$.
\end{remark}

\subsection{Boundary of the even Petrovskii class for $P_8$ singularities}
\label{obstrp8}

\begin{proposition}
For any singularity of the class $P_8^1$ or $P_8^2$, presented by a homogeneous
function $f(x,y,z)$ of degree 3, and for any its non-discriminant perturbation
$f_\lambda$, the boundary of the {\em even} local Petrovskii class is
non-trivial in $\tilde H_1(\partial V_\lambda)$.
\end{proposition}

{\it Proof.} The isomorphism class of the intersection form of the
singularities $P_8$ is well-known, see e.g. \cite{Ga1}. The group $\tilde
H_1(\partial V_\lambda)$ for any non-discriminant (and hence for any at all)
perturbation $f_\lambda$ of this function $f$ can be easily calculated from
(\ref{exact}) and is equal to $\Z^2 \oplus \Z_3$.

We can assume that the coordinates $x,y,z$ are chosen to take $f$ in the
Newton-Weierstrass normal form $x^3+axz^2 + bz^3 - y^2z$. The Hopf bundle
projection $S^5 \to \CP^2$ maps $\partial V_0$ to the elliptic curve $\{f=0\}$.
Regarding the submanifolds $f^{-1}(\tau) \cap B_\varepsilon \cap \R^n$, $0<
\tau <<\varepsilon$, realizing the classes $P_{\mbox{ev}}(f-\tau)$, and tending
$\tau$ to 0, we see that the class of $\partial P_{\mbox{ev}}(0)$ is mapped by
this Hopf projection into twice the class of the real part of this elliptic
curve, oriented as the boundary of the part of the affine chart $\{z=1\}$ in
$\RP^2$, in which the function $f$ takes negative values. The latter class is
never homological to zero in the elliptic curve. (In the case of $P_8^2$, when
this real part consists of two components, some sum of these components is
homological to zero, but the orientations of these components in this sum
should be coordinated in a different way.) \hfill $\Box$

\subsection{Boundary of $P_{\mbox{ev}}$ for critical points of functions
of two variables}

Denote by $D_\varepsilon$ the real part $B_\varepsilon \cap \R^n$ of the ball
$B_\varepsilon$. The real zero set of a function $f(x_1,x_2)$ with a critical
point at 0 consists of several irreducible curves passing through $0$. Any such
curve intersects the circle $\partial D_\varepsilon$ at two points. The
corresponding {\em chord diagram} is the graph consisting of the circle
$\partial D_\varepsilon$ and all its chords connecting the endpoints of any
such component. \label{obstr2dim}

\begin{proposition}
The boundary $\partial P_{\mbox{ev}} \in \tilde H_0(V_\lambda)$ is trivial for
some $($and then for all$)$ non-discriminant perturbations $f_\lambda$ of $f$
if and only if any chord of this chord diagram intersects an even number of
other chords.
\end{proposition}

{\em Proof.} Let us take the function $f-\tau$, $0<\tau <<\varepsilon$, for the
perturbation $f_\lambda$, so that its Milnor fiber is the set $f^{-1}(\tau)
\cap B_\varepsilon$. The geometric boundary of the set of real points of this
fiber is in the obvious one-to-one correspondence with the set of endpoints of
chords of the chord diagram. The points of this boundary belong to one and the
same component of the manifold $\partial V_\lambda$ if and only if they
correspond to the endpoints of one and the same chord. It is easy to calculate
that two points corresponding to the endpoints of some chord are counted for in
the homological boundary of the even Petrovskii cycle with one and the same
sign if and only if they are separated by an odd number of other endpoints in
the circle $\partial D_\varepsilon$. \hfill $\Box$

\section{Invariants of components of the complement of the real discriminant}
\label{maininv}

Let us choose the number $\Delta>0$ small enough so that all varieties $
f^{-1}(t),$ $t \in [-\Delta, \Delta]$, are transversal to $\partial
D_\varepsilon$.

Let $\Lambda \subset \R^l$ be a very small neighborhood of the origin in the
space of parameters $\lambda$, such that the same transversality condition is
satisfied not only for $f$, but also for all functions $f_\lambda$, $\lambda
\in \Lambda$, and additionally all real critical values of these functions
$f_\lambda$ in $D_\varepsilon$ belong to the interval $(-\Delta, \Delta)$.
Denote by $M_-(\lambda), M_0(\lambda)$ and $M_+(\lambda)$ the sets of lower
values $f_\lambda^{-1}((\infty,-\Delta]) \cap D_\varepsilon$,
$f_\lambda^{-1}((\infty,0]) \cap D_\varepsilon$, and
$f_\lambda^{-1}((\infty,\Delta]) \cap D_\varepsilon$, respectively.

The diagrams of spaces
\begin{equation}
\label{diag1}
\begin{array}{ccccc}
M_-(\lambda) & \subset & M_+(\lambda) & \subset & D_\varepsilon \\
\cup & & \cup & & \cup \\
M_-(\lambda) \cap \partial D_\varepsilon& \subset & M_+(\lambda) \cap
\partial D_\varepsilon & \subset & \partial D_\varepsilon
\end{array}
\end{equation}
form a locally trivial (and hence trivializable) fiber bundle over the
neighborhood $\Lambda$. Therefore we can fix a family (depending continuously
on $\lambda$) of homeomorphisms of all of them to one and the same diagram
corresponding to some distinguished value $\lambda_0$ of $\lambda$, say to
$\lambda_0=0$. The spaces $M_-(\lambda_0)$ and $M_+(\lambda_0)$ for this
distinguished value will be called just $M_-$ and $M_+$. Given an arbitrary
$\lambda$, composing the embedding $M_0(\lambda) \to D_\varepsilon$ with this
unifying homeomorphism we obtain the diagram of spaces

\begin{equation}
\label{diag12}
\begin{array}{ccccccc}
M_- & \subset & M_0(\lambda) & \subset &
M_+ & \subset & D_\varepsilon \\
\cup & & \cup & & \cup & & \cup \\
M_- \cap \partial D_\varepsilon & \subset & M_0(\lambda) \cap
\partial D_\varepsilon & \subset & M_+ \cap \partial D_\varepsilon &
\subset & \partial D_\varepsilon
\end{array}
\end{equation}

\begin{proposition}
\label{homotinv} If two points $\lambda, \lambda' \in \Lambda$ belong to one
and the same connected component of the set of non-discriminant perturbations
of $f$, then the corresponding diagrams $($\ref{diag12}$)$ are isotopic to one
another via an isotopy of the pair $(D_\varepsilon, \partial D_\varepsilon)$
constant on $M_-$ and on $D_\varepsilon \setminus M_+$. \hfill $\Box$
\end{proposition}

In particular, all homological invariants of isotopy classes of such diagrams
also are invariants of components of the complement of the discriminant, and we
get the following corollary.

\begin{proposition} \label{homolinv}
The following objects are the same for all $\lambda$ from one and the same
component of the complement of the discriminant:

a$)$ the isomorphism classes of groups $H_*(M_0(\lambda))$, $H_*(M_0(\lambda),
\partial D_\varepsilon)$, $H_*(M_0(\lambda), M_-(\lambda))$,
$H_*(M_0(\lambda),(M_- \cup \partial D_\varepsilon))$, $H_*(M_+,M_0(\lambda))$,
and $H_*(M_+,(M_0(\lambda)\cup \partial D_\varepsilon))$;

b$)$ images of boundary operators $\partial: H_*(M_0(\lambda),M_-) \to
H_*(M_-)$ and $\partial: H_*(M_0(\lambda), M_- \cup
\partial D_\varepsilon) \to H_*(M_- \cup \partial D_\varepsilon)$,

c$)$ kernels of operators defined by inclusions, $H_*(M_-) \to
H_*(M_0(\lambda))$, $H_*(M_- \cup \partial D_\varepsilon) \to H_*(M_0(\lambda)
\cup
\partial D_\varepsilon)$, $H_*(M_+, M_-) \to H_*(M_+, M_0(\lambda))$, etc.
\hfill $\Box$
\end{proposition}

The invariant $\chi(\lambda)$ used in Proposition \ref{chi} is just the Euler
characteristic of the third group mentioned in item (a) of Proposition
\ref{homolinv}.

\section{Realization of local lacunas promised in Theorem \ref{thm1}}
\label{realiz}

\subsection{Classes $P_8^1$ and $P_8^2$} \label{realp8}

\begin{proposition}
If $f(x_1, x_2, x_3)$ is a non-degenerate homogeneous polynomial of degree 3
$($so that it belongs to one of classes $P_8^1$ or $P_8^2)$ then the
polynomials $f_{\pm \tau} \equiv f \pm (\tau (x_1^2+x_2^2+x_3^2) - \tau^3)$
with sufficiently small $\tau>0$ belong to odd local lacunas.

Moreover, the perturbations $f_\tau$ and $f_{-\tau}$ belong to {\em different}
odd local lacunas.
\end{proposition}

{\it Proof.} The odd Petrovskii cycle of $f_{\pm \tau}$, realized as in \S
\ref{another} with the central point ${\bf x}$ at the origin, is empty. Indeed,
any complex line through $0$, which is the complexification of a real line,
intersects $V_{\pm \tau}$ in at least two real points. The set of non-real
intersection points is complex conjugate to itself and consists of no more than
one point, since the degree of $f_{\pm \tau}$ is equal to 3.

The perturbations $f_\tau$ and $f_{-\tau}$ are separated by an invariant from
Proposition \ref{homolinv}(a). Namely, the relative homology group
$H_*(M_0(\tau), (M_- \cup
\partial D_\varepsilon))$ coincides with the homology group of a single point,
and the group $H_*(M_0(-\tau), (M_- \cup \partial D_\varepsilon))$ is
isomorphic to $H_*(S^2,\mbox{pt})$. \hfill $\Box$

\subsection{The class $\pm X_9$}
\label{lex}

We will consider the singularity class $+X_9$ only, since the class $-X_9$ can
be reduced to it, see \S \ref{pm}.

The local lacuna for a singularity of the class $+X_9$ assumed in the second
column of Table \ref{t13} is described in \S \ref{examin}. One of two lacunas
assumed in the fourth column is described in Corollary \ref{extra} on page
\pageref{extra} and is represented by the function $\varphi(x_1, x_2) -
x_3^2+\tau$, where $\tau>0$ is small enough and $\varphi$ has a minimum point
at $0$.

\begin{proposition}
\label{otherx9} If the function $\varphi(x_1,x_2)$ of the class $+X_9$ is a
non-negative homogeneous polynomial of degree 4, then the function $f_\tau
\equiv \varphi(x_1,x_2)-\tau(x_1^2+x_2^2)-x_3^2+ \tau^3$ with sufficiently
small $\tau>0$ belongs to the odd local lacuna of the function
$\varphi(x_1,x_2)-x_3^2$. This lacuna is different from the second lacuna
indicated in the previous paragraph.
\end{proposition}

{\it Proof.} By the Morse lemma, changing slightly the local coordinate $x_3$
(which certainly does not change the values of the Petrovskii classes) we can
replace the function $-x_3^2$ in one variable by $-x_3^2+x_3^4$, and hence the
function $f_\tau$ by $f_\tau+x_3^4$. The corresponding odd Petrovskii cycle,
described in \S \ref{another} for ${\bf x}=0$, is empty since any real line
through $0$ intersects the zero set of this function $f_\tau +x_3^4$ at four
points. The last statement of the proposition follows immediately from
Proposition \ref{homolinv}. \hfill $\Box$

\subsection{Remaining lacunas for corank 2 parabolic singularities}

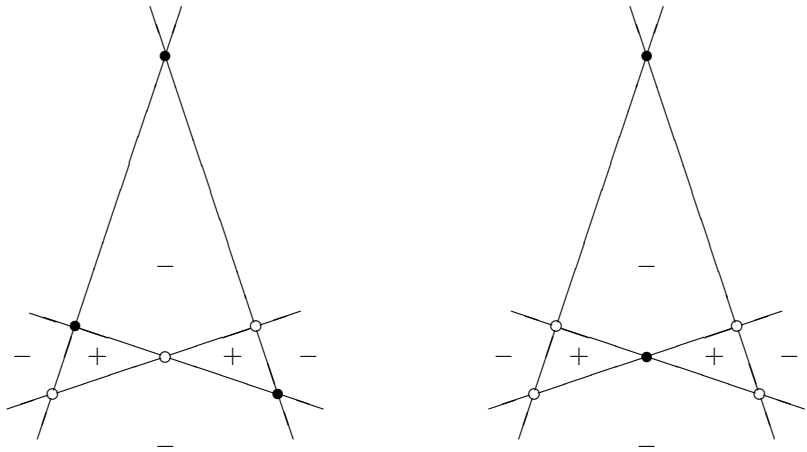
\begin{figure}
\unitlength 1.00mm \linethickness{0.6pt}
\begin{center}
\begin{picture}(106.00,65.00)
\put(2,18){\makebox(0,0)[cc]{$-$}} \put(40,18){\makebox(0,0)[cc]{$-$}}
\put(21,6){\makebox(0,0)[cc]{$-$}}

\put(66,18){\makebox(0,0)[cc]{$-$}} \put(104,18){\makebox(0,0)[cc]{$-$}}
\put(85,6){\makebox(0,0)[cc]{$-$}}

\put(21.00,18.00){\circle{1.33}} \put(33.00,22.00){\circle{1.33}}
\put(36.00,13.00){\circle*{1.33}} \put(6.00,13.00){\circle{1.33}}
\put(9.00,22.00){\circle*{1.33}} \put(21.00,58.00){\circle*{1.33}}
\put(0,11){\line(3,1){5.40}} \put(6.60,13.20){\line(3,1){13.80}}
\put(21.60,18.20){\line(3,1){10.80}} \put(33.60,22.20){\line(3,1){5.40}}
\put(42.00,11.00){\line(-3,1){5.40}} \put(35.40,13.20){\line(-3,1){13.80}}
\put(20.40,18.20){\line(-3,1){10.80}} \put(8.40,22){\line(-3,1){5.40}}
\put(4.00,7.00){\line(1,3){1.80}} \put(6.20,13.60){\line(1,3){2.60}}
\put(9.20,22.60){\line(1,3){11.60}} \put(21.20,58.60){\line(1,3){2}}
\put(19.00,64.60){\line(1,-3){2}} \put(21.20,57.40){\line(1,-3){11.60}}
\put(33.20,21.40){\line(1,-3){2.60}} \put(36.20,12.40){\line(1,-3){1.80}}
\put(21.00,30.00){\makebox(0,0)[cc]{$-$}}
\put(30.00,18.00){\makebox(0,0)[cc]{$+$}}
\put(12.00,18.00){\makebox(0,0)[cc]{$+$}} \put(85.00,18.00){\circle*{1.33}}
\put(97.00,22.00){\circle{1.33}} \put(100.00,13.00){\circle{1.33}}
\put(70.00,13.00){\circle{1.33}} \put(73.00,22.00){\circle{1.33}}
\put(85.00,58.00){\circle*{1.33}} \put(64.00,11.00){\line(3,1){5.40}}
\put(70.60,13.20){\line(3,1){13.80}} \put(85.60,18.20){\line(3,1){10.80}}
\put(97.60,22.20){\line(3,1){5.40}} \put(106.00,11.00){\line(-3,1){5.40}}
\put(99.44,13.20){\line(-3,1){13.80}} \put(84.40,18.20){\line(-3,1){10.80}}
\put(72.40,22.20){\line(-3,1){5.40}} \put(68.00,7.00){\line(1,3){1.80}}
\put(70.20,13.60){\line(1,3){2.60}} \put(73.20,22.60){\line(1,3){11.60}}
\put(85.20,58.60){\line(1,3){2}} \put(82.80,64.60){\line(1,-3){2}}
\put(85.20,57.40){\line(1,-3){11.60}} \put(97.20,21.40){\line(1,-3){2.60}}
\put(100.20,12.40){\line(1,-3){1.80}} \put(85.00,30.00){\makebox(0,0)[cc]{$-$}}
\put(94.00,18.00){\makebox(0,0)[cc]{$+$}}
\put(76.00,18.00){\makebox(0,0)[cc]{$+$}}
\end{picture}
\caption{Lacunas for $X_9^2$} \label{X92}
\end{center}
\end{figure}

By Proposition \ref{stabil} all remaining local lacunas assumed in non-zero
cells of Table \ref{t13} can be considered as odd local lacunas of some
functions in two variables. We realize these lacunas in the following way. As
in \cite{Gusein-Zade 74}, \cite{A'Campo 75}, we demonstrate a perturbation
$f_\lambda(x_1,x_2)$ of the corresponding function $f$, all whose $\mu(f)$
critical points are real, all critical values at the saddlepoints are equal to
0, and all critical values at minima (respectively, maxima) are negative
(respectively, positive). Using a further very small perturbation  of
$f_\lambda$ we can obtain a function $f_{\tilde \lambda}$ arbitrarily close to
$f_\lambda$ but with critical values at all saddlepoints moved from 0 to any
prescribed sides; in particular $f_{\tilde \lambda}$ is non-discriminant. In
the Figures \ref{X92}, \ref{j103} we draw the zero sets of the preliminary
perturbations $f_\lambda$, and indicate by black (respectively, white) circles
the saddlepoints, the values at which should be moved in the negative
(respectively, positive) direction from 0.

\unitlength 0.80mm
\begin{figure}
\begin{center}
\begin{picture}(100,47)
\put(0,20){\line(1,0){63}} \put(83,20){\line(-1,0){18}}
\put(85,20){\line(1,0){9}} \bezier{120}(5,0)(15,24)(27.5,31.4)
\bezier{30}(29.5,32.2)(34,35)(40,35) \bezier{100}(40,35)(53,35)(63.1,21)
\bezier{50}(65,19)(71,10)(75,0) \bezier{88}(60,5)(73,5)(83,19)
\bezier{50}(84.6,21)(90,30)(95,40) \bezier{120}(60,5)(40,5)(29,31)
\bezier{24}(28.2,33)(26.3,37)(25.6,40) \put(28.70,32){\circle{1.7}}
\put(64,20){\circle{1.7}} \put(84,20){\circle{1.7}} \put(16,20){\circle*{1.7}}
\put(35,20){\circle*{1.7}} \put(71.4,8){\circle*{1.7}} \put(24,24){{\Large
$-$}} \put(42,25){{\Large $+$}} \put(54,11){{\Large $-$}} \put(71,13){{\Large
$+$}}
\end{picture}
\end{center}
\caption{Lacuna for $J_{10}^3$} \label{j103}
\end{figure}
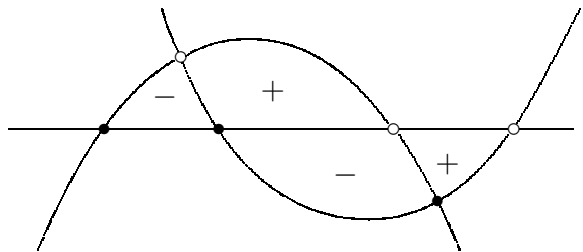

\begin{proposition}
If a function $f(x_1,x_2)$ of the class $X_9^2$ is represented by a homogenous
polynomial of degree 4 vanishing on four different real lines, then

1$)$ it has a perturbation $f_\lambda$ whose zero set is as shown in either
side of Fig.~\ref{X92};

2$)$ the further non-discriminant perturbations $f_{\tilde \lambda}(x_1,x_2)$
shown in these pictures by black and white circles belong to odd local lacunas
of $f$;

3$)$ these two local lacunas are different;

4$)$ rotating both pictures of Fig.~\ref{X92} by the angle $\pi/2$ we obtain
the pictures of two other perturbations of $f$ which belong to two additional
odd local lacunas different from the previous two.
\end{proposition}

{\it Proof.} Statement 1 is obvious, 2 follows from the calculation of odd
Petrovskii cycles mentioned in \S \ref{expl}, and statements 3, 4 follow from
Proposition \ref{homolinv}: indeed, the images of the boundary operators
$H_1(M_0(\lambda),M_-) \to H_0(M_-) \simeq \Z^4$ for these four cases are four
different subgroups of the latter group. \hfill $\Box$

\begin{proposition}
If the function $f(x_1,x_2)$ belongs to the class $J_{10}^3$, then

1$)$ it has a perturbation $f_\lambda$ whose zero set is homeomorphic to the
one shown in Fig.~\ref{j103},

2$)$ the further non-discriminant perturbation $f_{\tilde \lambda}(x_1,x_2)$
described in this picture by black and white circles belongs to an odd local
lacuna of $f$.
\end{proposition}

{\it Proof.} This proposition follows immediately from the normal form of
critical points of type $J_{10}^3$ and from the calculation of odd Petrovskii
cycles discussed in \S \ref{expl}. \hfill $\Box$

\section{A program counting for topologically different
morsifications of critical points of real functions} \label{progr}

This program has two versions: one for singularities of corank $\leq 2$ (see \\
\verb"https://www.hse.ru/mirror/pubs/share/185895886", currently it contains
the starting data for the singularity class $J_{10}^3$), and the other one for
singularities of arbitrary ranks
(\verb"https://www.hse.ru/mirror/pubs/share/185895827", currently with initial
data of $P_8^1$). The further versions of the program will occur at the bottom
of the page  \verb"https://www.hse.ru/en/org/persons/1297545#sci".

For a description of the program see \S V.8 of the book \cite{APLt}, however
the web reference given there leads to an obsolete version of the program.

The starting data for the program are the topological characteristics of some
morsification $f_\lambda$ of $f$, all whose critical points are real, and all
their critical values are different and not equal to $0$. Namely, these data
include the Morse indices of all critical points ordered by the increase of
their critical values (in the program for corank=2 singularities), or just the
parities of these indices (in the program for for the general case) and the
intersection indices of corresponding vanishing cycles in $\tilde
H_{n-1}(V_\lambda)$, defined by a canonical system of paths and having
canonical orientations compatible with the orientation of $\R^n$. One
additional element of data is the number of negative critical values of
$f_\lambda$. This information is sufficient to calculate both Petrovskii
classes of our morsification $f_\lambda$ and of all its stabilizations.

Our program is modelling (on the level of similar sets of topological data) all
potentially possible topological surgeries of the initial morsifications,
namely the jumps of critical values through $0$, collisions of the real
critical values (which can either bypass one another or undergo a Morse surgery
and go into the imaginary domain), the opposite operations (that is, a
collision of two complex conjugate critical values at a real point), and also
rotations of imaginary critical values around one another. Knowing our
topological data before any of these surgeries is enough to predict the similar
data after it.

In general, it is not sure that any sequence of such operations over the sets
of topological invariants actually can be realized by a path in the parameter
space of the deformation, so we consider their results as {\em virtual
morsifications}, that is, some admissible collections of our topological data,
including the Petrovskii classes. However, any actual morsification is
represented by a virtual one, which surely will be found by our algorithm (if
it will have enough of memory and time). In particular, if the program has
enumerated all possible virtual morsifications of a singularity class and found
that their Petrovskii classes never vanish, then we can put the sign $0_c$ in
the corresponding cell of Table \ref{t13}.

On the other hand, a majority of real local lacunas described in \S
\ref{realiz} was discovered by this program. More precisely, it has found
suspicious virtual morsifications with vanishing Petrovskii classes and printed
out their topological data; after that in all our cases it was easy to find by
hands the real morsifications with these data.

The numbers of topologically distinct virtual morsifications found by our
program are equal to 6503 for $P_8^1$, 9174 for $P_8^2$, 16928 for $\pm X_9$,
96960 for $X_9^2$, 549797 for $J_{10}^1$, and 77380 for $J_{10}^3$.

\section{Starting data for the singularity $P_8^2$}
\label{DDP82}

The initial data of our program (that is, convenient morsifications with only
real critical points and intersection indices of their vanishing cycles) for
all real parabolic singularities except for $P_8^2$ can be easily calculated by
the methods of \cite{Gusein-Zade 74}, \cite{A'Campo 75} (for critical points
whose quadratic part is of corank $\leq 2$) or \cite{Ga1} (for the class
$P_8^1$ which has a convenient representative $x^3+y^3+z^3$). It is important
for our algorithm that the orientations of these vanishing cycles, defining the
signs of the intersection indices, should be compatible with the fixed
orientation of $\R^n$; fortunately all these methods satisfy this condition. In
this section we solve the similar problem for the remaining case $P_8^2$.

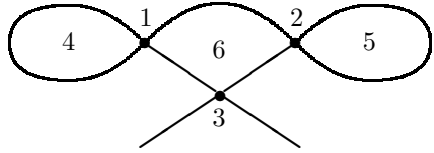
\begin{figure}
 \unitlength 1.00mm
\linethickness{0.4pt} \thicklines
\begin{picture}(56.00,27.00)
\put(17.3,19.3){$1$} \put(37.3,19.3){$2$} \put(27,6){$3$} \put(7,16){$4$}
\put(47,16){$5$} \put(27,15){$6$}

\put(18.00,17.00){\circle*{1.33}} \put(38.00,17.00){\circle*{1.33}}
\bezier{100}(18.50,17.50)(28.00,27.00)(37.50,17.50)
\put(28.00,10.00){\circle*{1.33}}

\bezier{44}(17.50,16.50)(13.00,12.00)(8.00,12.00)
\bezier{44}(17.50,17.50)(13.00,22.00)(8.00,22.00)
\bezier{52}(8.00,12.00)(0.00,12.00)(0.00,17.00)
\bezier{52}(0.00,17.00)(0.00,22.00)(8.00,22.00)

\bezier{44}(38.50,16.50)(43.00,12.00)(48.00,12.00)
\bezier{44}(38.50,17.50)(43.00,22.00)(48.00,22.00)
\bezier{52}(48.00,12.00)(56.00,12.00)(56.00,17.00)
\bezier{52}(56.00,17.00)(56.00,22.00)(48.00,22.00)

\put(27.40,9.80){\line(-3,-2){10}} \put(28.60,9.80){\line(3,-2){10}}
\put(27.40,10.60){\line(-3,2){9.20}} \put(28.60,10.60){\line(3,2){9.20}}
\end{picture}
\caption{Convenient morsification for $E_6$} \label{sabe6}
\end{figure}

We choose a function in this class whose Newton--Weierstrass normal form is $f
= x^3-x z^2 + y^2z$. Consider its small perturbation $f_1=f+ \varepsilon z^2$,
$\varepsilon>0$; by the dilation of the coordinates and the function we can
assume that $\varepsilon=1$. $f_1$ has a critical point of type $E_6$ at the
origin. By a small (with unit linear part) local change of coordinates at the
origin this function can be reduced to the form $\tilde x^3 - y^4 + \tilde
z^2$. In addition $f_1$ has two real Morse critical points $(1,0,\sqrt{3})$ and
$(1,0,\sqrt{3})$ with common critical value $1$; their Morse indices are equal
to 2 and 1 respectively. As was shown in \cite{Gusein-Zade 74}, Example 3, we
can slightly perturb our function $f_1$ so that its $E_6$-type critical point
splits into six real Morse critical points, and the new function $f_2$ has the
form $\varphi(\tilde x,\tilde y) + \tilde z^2$ in the corresponding
neighborhood of the origin, where the zero level set of $\varphi$ in $\R^2$
looks as in Fig.~\ref{sabe6}. The crossing points of this set correspond to the
Morse critical points of $f_2$ with critical value 0 and Morse index 1, and any
of three bounded domains contains a point with a slightly greater critical
value and Morse index 2. Let $f_3$ be an additional very small perturbation of
the Morse function $f_2$ making it a strictly Morse function, that is,
separating all 8 critical values. We can do it in such a way that the order of
these values in $\R^1$ will be as indicated by the numbers in Fig.~\ref{sabe6}.
Choose a real value $A$ greater than all six critical values obtained from the
perturbation of $E_6$-type critical point but lower than the critical values at
two Morse critical points obtained from two points with critical value 1.
Define the basis of vanishing cycles in $H_2(f_3^{-1}(A))$ by the system of
paths connecting $A$ to all eight critical values within the upper half-plane
in $\C^1$. Number the first six basis cycles $\Delta_i$ by the order of
corresponding critical values, see Fig. \ref{sabe6}; let the 7th and the 8th
cycles be the ones arising from the critical points with value $\approx 1$ and
Morse indices 2 and 1 respectively. Orient all these cycles in the
correspondence with the fixed orientation of $\R^n$ (see page 177 in
\cite{APLt}, especially formula (V.6) there). Our purpose is to calculate the
matrix of intersection indices of these eight cycles: this will give us a set
of initial data for the perturbation $f_3-A$ of the initial function $f$. This
matrix is symmetric since $n-1$ is even.

\begin{equation} \left(
\begin{array}{cccccccc}
-2 \ & 0 & 0 & 1 & 0 & 1 & X & -\frac{Z+W}{2} \\
0 & -2 \ & 0 & 0 & 1 & 1 & X & -\frac{Z+W}{2} \\
0 & 0 & -2 \ & 0 & 0 & 1 & Y & -\frac{W}{2} \\
1 & 0 & 0 & -2 \ & 0 & 0 & 0 & Z \\
0 & 1 & 0 & 0 & -2 \ & 0 & 0 & Z \\
1 & 1 & 1 & 0 & 0 & -2 \ & 0 & W \\
X & X & Y & 0 & 0 & 0 & -2 \ & 0 \\
-\frac{Z+W}{2} & -\frac{Z+W}{2} & -\frac{W}{2} & Z & Z & W & 0 & -2 \
\end{array}
\right) \label{matp8}
\end{equation}

\begin{lemma}
The wanted intersection matrix has the form $($\ref{matp8}$)$ for some values
$X, Y, Z$ and $W$. \end{lemma}

{\it Proof.} The intersection indices of the first six cycles can be calculated
by the method of \cite{Gusein-Zade 74}, \cite{A'Campo 75} and are as shown in
the upper left-hand $6 \times 6$ corner of the matrix (\ref{matp8}). The
intersection index $\langle \Delta_7, \Delta_8\rangle$ is equal to 0 because
these cycles appear from distant critical points with almost coinciding
critical values. The cycles $\Delta_4, \Delta_5$ and $\Delta_6$ are invariant
under the complex conjugation in $f_3^{-1}(A)$, therefore their intersection
indices with $\Delta_7$ (which is anti-invariant) are equal to 0, as indicated
in the 7th row of (\ref{matp8}). Also, the perturbation $f_1$ of the original
function $f$ is invariant under the reflection in the hyperplane $y=0$, and its
further perturbation $f_2$ can be accomplished keeping this symmetry. This
reflection keeps the basis cycles in $f_2^{-1}(A)$ which are close to the
cycles $\Delta_7$ and $\Delta_8$, only changing their canonical orientations;
on the other hand it permutes the cycles close to $\Delta_1$ and $-\Delta_2$,
and also cycles close to $\Delta_4$ and $-\Delta_5$. Therefore $\langle
\Delta_1,\Delta_7\rangle =\langle-\Delta_2,-\Delta_7\rangle \equiv
\langle\Delta_2,\Delta_7\rangle$ and $\langle \Delta_4,\Delta_8\rangle
=\langle\Delta_5,\Delta_8\rangle$. It is why the corresponding cells of our
matrix (\ref{matp8}) are filled in by equal letters ($X$ in the first case and
$Z$ in the second).

For any $i=1,2,3$ denote by $\bar \Delta_i$ the vanishing cycle in
$f_3^{-1}(A)$ obtained from the $i$th critical point by the path connecting the
corresponding critical value with $A$ in the {\em lower} half-plane of $\C^1$.
It is easy to see that the cycle $\Delta_i+ \bar \Delta_i$ is anti-invariant
under the complex conjugation and $\Delta_8$ is invariant, therefore $\langle
\Delta_i + \bar \Delta_i, \Delta_8\rangle =0$. But $\bar \Delta_i$ can be
considered as the image of $\Delta_i$ under the Picard-Lefschetz monodromy
operator along a loop $L$ starting and ending at the point $A$ and embracing
all critical values placed between the $i$th one and the point $A$. This image
is equal to $\Delta_i + \mbox{Var}_L(\Delta_i)$, therefore the previous
equation gives us $-2\langle \Delta_i,\Delta_8\rangle = \langle
\mbox{Var}_L(\Delta_i),\Delta_8 \rangle$. By the Picard--Lefschetz formula the
last number is equal to $\sum_{j=4}^6 \langle \Delta_i,\Delta_j\rangle \langle
\Delta_j, \Delta_8 \rangle$, which gives us the expressions of intersection
indices $\langle \Delta_i,\Delta_8\rangle$ through the numbers
$Z\equiv\langle\Delta_{4},\Delta_8\rangle = \langle\Delta_{5},\Delta_8\rangle$
and $W\equiv\langle\Delta_6,\Delta_8\rangle$; see the first three cells of the
last row of (\ref{matp8}). \hfill $\Box$ \medskip

It remains to calculate the numbers $X, Y, Z$ and $W$ in this matrix. We know
that it is the matrix of the $P_8$ type bilinear form in some basis in $\Z^8$.
This form is well-known, see e.g. \cite{Ga1}. In particular it is easy to check
that the image of this lattice in the dual lattice under the map defined by
this bilinear form coincides with the image of an arbitrary $E_6$-sublattice.
Therefore the last two rows of the desired matrix are integer linear
combinations of the first six ones. Writing these rows in the form of such
linear combinations with indeterminate coefficients, we get 16 equations in 16
unknowns $a_1, \dots, a_6, b_1, \dots, b_6, X, Y, Z, W$. Some two of these
equations are consequences of the others, but the diophantine system of
remaining 14 equations is easily solvable and has exactly four different
integer solutions, which imply four possible combinations of coefficients of
the matrix (\ref{matp8}): $\{X=0, Y=\pm 1$, $Z=0, W = \pm 2\}$. Let us select
the correct version.

Both local Petrovskii classes of a morsification, all whose critical points are
real, can be calculated explicitly from the Morse indices of these critical
points and intersection indices of all (properly oriented) vanishing cycle, see
\S \ref{expl}. Substituting all four hypothetical combinations of intersection
indices in these calculations, in three cases we get a contradiction with the
previously obtained results on these classes for $P_8^2$ singularities, saying
that

a) in the case of odd $n$ and even index $i_+$ (e.g. for $n=3$) these
singularities have local lacunas, therefore the homological boundary of the odd
Petrovskii class is equal to 0 for all non-discriminant perturbations of $f$,
see \S \ref{realp8};

b) the similar homological boundary of the even Petrovskii class is not equal
to 0 for the same values of $n$ and $i_+$, see \S \ref{obstrp8}.

The unique remaining case gives us $Y=1$, $W=-2$, and the intersection matrix
is completely calculated. Plugging it into the initial data of our program, we
get from it the messages that no local lacunas exist close to the $P_8^2$
singularities in the cases of even $n$ and both even and odd $i_+$. This proves
the first two zeros in Table \ref{t13} for $P_8^2$.

\section{No extra lacunas for $\pm X_9$}
\label{noextra}

The sign $1_{c}$ in the second column of Table \ref{t13} for $\pm X_9$ (that
is, the fact that this singularity has no local lacunas in addition to the one
mentioned in \S \ref{lex}) is proved by our program with the help of the
following fact.

\begin{proposition}
If the function $f$ has a minimum point at the origin, then all its
sufficiently small perturbations $f_\lambda$, such that all real critical
values of $f_\lambda$ are positive, belong to one and the same component of the
complement of the discriminant of an arbitrary its versal deformation.
\end{proposition}

{\it Proof.} The described property of functions is preserved by inducing and
equivalence of deformations, therefore it is enough to prove our proposition
for one arbitrary versal deformation of the function $f$. In particular we can
assume that this deformation contains together with any perturbation
$f_\lambda$ of $f$ also all perturbations $f_\lambda + c$, where the constants
$c$ run some interval containing 0. Choose some small value $c>0$ in this
interval, then there is some number $\delta>0$ such that the
$\delta$-neighborhood of the point $\{f+c\}$ in the space $\R^l$ of parameters
of this deformation is separated from the discriminant. For any point $\lambda$
from the $\delta$-neighborhood of the origin in $\R^l$, such that $f_\lambda$
satisfies the condition of our Proposition, the entire segment consisting of
functions $f_\lambda+\tau$, $\tau \in [0,c]$ belongs to the complement of the
discriminant, and its last point $f_\lambda+c$ belongs to the
$\delta$-neighborhood of the point $f+c$ mentioned above. \hfill $\Box$
\medskip

Therefore I have asked my program to check that the functions of type $+X_9$ do
not have virtual morsifications with trivial even Petrovskii class and at least
one negative critical value. Its confirmative answer justifies the sign $1_c$
in the cell under question.

\end{document}